\documentclass[10 pt, reqno]{amsart}
\usepackage{amsmath}
\usepackage{amssymb,graphicx,color,amsthm,hyperref,subfigure,amsbsy}
\usepackage{subfigure}

\let\mib=\boldsymbol
\let\cal=\mathcal

\def\C{{\bf C}}
\def\N{{\bf N}}
\def\R{{\bf R}} 
\def\vu{{\sf u}}
\def\mx{{\bf x}} 
 
\def\mA{{\bf A}} 
\def\mI{{\bf I}} 
\def\mN{{\bf N}} 
\def\malpha{{\mib \alpha}}
\def\mmu{{\mib \mu}}
 
\def\mxi{{\mib \xi}}

\def\vu{{\sf u}} 
\def\vu{{\sf v}}

\def\eps{\varepsilon} 
\def\lA{{\cal A}}
\def\lF{{\cal F}}
\def\lS{{\cal S}}

\def\Bnor#1#2{\Bigl\| #1 \Bigr\|_{#2}}
\def\CC{{\rm C}}
\def\H#1#2{{{\rm H}^{#1}(#2)}}
\def\Hl#1#2{{{\rm H}^{#1}_{{\rm loc}}(#2)}} 
\def\Hnl#1#2{{{\rm H}^{#1}_{0}(#2)}} 
\def\HH#1{{{\rm H}^{#1}}}

\def\dv{{\sf div\thinspace}}
\def\dscon{\relbar\joinrel\rightharpoonup}
\def\dstr{\longrightarrow}
\def\Dup#1#2{\langle#1,#2\rangle}
\def\Ld#1{{{\rm L}^{2}(#1)}}
\def\Ldl#1{{{\rm L}^{2}_{{\rm loc}}(#1)}}
\def\LLd{{{\rm L}^2}}
\def\mnul{{\bf 0}} 
\def\nor#1#2{{\| #1 \|}_{#2}}

\def\oi#1#2{\langle#1,#2\rangle}
\def\ozi#1#2{\langle#1,#2]}
\def\paragraf{\mathhexbox278}
\def\Pd{{\rm P}}
\def\Rd{{{\bf R}^{d}}}
\def\Rjpd{{{\bf R}^{1+d}}}
\def\Rpl{{{\bf R}^{+}}} 
\def\Sd{{\rm S}^{d}}
\def\str{\longrightarrow}
\def\supp{{\rm supp\,}}
\def\vnul{{\sf 0}}

\newtheorem{theorem}{Theorem}[section]
\newtheorem{lemma}[theorem]{Lemma}
\newtheorem{proposition}[theorem]{Proposition}

\newtheorem{remark}{Remark}
\newtheorem{example}{Example}

\numberwithin{equation}{section} %za numeriranje jed-bi s brojem sekcije

\begin{document}

\title[Stability of Observations of PDE-s under Uncertain Perturbations]{Stability of Observations of Partial Differential Equations under Uncertain Perturbations}

\author{Martin Lazar}

\address{Martin Lazar, University of Dubrovnik, Department of Electrical Engineering and Computing,
\'Cira Cari\' ca 4, 20 000 Dub\-rov\-nik, Croatia \\
phone: ++385\,20\,445842, \quad fax: ++385\,20\,435590 }
\email{martin.lazar@unidu.hr} 
\urladdr{http://www.martin-lazar.from.hr}
%\phone{++385 20 445842} 
%\fax{++385 20 435590} 

\begin{abstract}
We analyse stability of observability estimates for solutions to wave and Scr\" odinger equations subjected to  additive perturbations. 
The paper  generalises the recent averaged observability/control result  by allowing for systems consisting of operators of different types. 
The method also applies to the simultaneous observability problem by which one tries to estimate the energy of each component of a system under consideration.
The analysis relies on microlocal defect tools; in particular on standard H-measures, when the main dynamic of the system is governed by the wave operator, while parabolic H-measures are explored in the case of the Schr\" odinger operator.

%The method of proof relies either on standard or parabolic H-measures, in dependence on the operator, the wave or the Schr\" odinger one,  governing  the main dynamic of the system under considerations.
\end{abstract}

\subjclass{93B05, 93B07,  93C20, 93D09 }

\keywords{averaged control, robust observability, parabolic H-measures}

\maketitle

\section{Introduction}

A notion of averaged control has been recently introduced in \cite{Z-aver, LZ}, both for parameter dependent ODEs, as well as for systems of PDEs with variable coefficients. Its goal is to control the average (or more generally a suitable linear combination) of system components by a single control. The problem is relevant in practice where the control has to be chosen independently of the coefficient value.

The notion is equivalent to the averaged observability, by which the energy  of the system is recorded by observing the average of solutions on a suitable subdomain. 

In the paper we investigate a more general problem based on a system whose first component represents its main dynamic, while the other ones correspond to perturbations. Assuming that the main component is observable, we explore conditions by which that property remains stable under additive perturbations.

In general, the operators entering the system are not assumed to be of the same type. In a special case of a system consisting of a same type operators, the result corresponds to the averaged control of the system, thus incorporating results obtained in \cite{LZ}. 

The methods are applied to the simultaneous observability problem as well, by which one tries to estimate energy of all system components by observing their average. The  corresponds dual problem consists of controlling each individual component of the adjoint system by means of a same control.

The study of the problem explores microlocal analysis tools, in particular H-measures and their variants. H-measures, introduced independently in \cite{Ger,Tar}, are kind of defect tools, measuring deflection of  the weak from strong convergence of $\LLd$ sequences. Since their introduction, they have been successfully applied in many mathematical fields - let us just mention generalisation of compensated compactness results to equations with variable coefficients \cite{Ger,Tar} and applications in the control theory \cite{BG97, DLR, LZ}.  Most of these applications apply  the so called {\it localisation principle} providing constraints on the support of H-measures (e.g. \cite{Tar}), and the proofs of this paper rely on it as well.

The paper is organised as follows. In the next section  we provide an averaged observability result for  a system whose main dynamic is governed by the wave operator. The finite system is analysed first, followed by  generalisations to an infinite discrete setting. Application of the approach to simultaneous observability is provided in the subsection 2.3.
The third section is devoted to observation of the Schr\"odinger equation under perturbations determined either by a hyperbolic or by a parabolic type operator. In the latter case, parabolic H-measures (generalisation of original ones to a parabolic setting) have to be explored. 
The paper is closed with concluding remarks, and by pointing toward some open and related problems.

\section{Observation of the wave equation under uncertain perturbations}
\label{wave}
\subsection{Averaged observability}

We analyse the problem of recovering the energy of the wave equation by observing an additive perturbation of the solution. The perturbation 
is  determined by a differential operator $P_2$, in general different from the wave one.

More precisely, we consider  the following system of equations:
\begin{equation}
\label{discrete-wave}
\begin{split}
P_1 u_1= \partial_{tt} u_1 - \dv (\mA_1(t, \mx) \nabla u_1)&=0, \qquad (t, \mx) \in \Rpl\times\Omega\cr
P_2 u_2&=0, \qquad (t, \mx) \in \Rpl\times\Omega\cr
u_1&=0, \qquad (t, \mx) \in \Rpl\times\partial\Omega\cr
u_1(0, \cdot)&=\beta_0 \in \Ld\Omega\cr
  \partial_t u_1(0, \cdot)&=\beta_1\in \H{-1}\Omega \,,\cr
\end{split}
\end{equation}
where  $\Omega$ is an open, bounded set in $\Rd$, $\mA_1$ is a bounded, positive definite matrix field, while $P_2$ is some, almost arbitrary, differential operator (precise conditions on it will be given below).  In the sequel we shall also use the notation $\lA_1= - \dv (\mA_1 \nabla)$ for the elliptic part of $P_1$. 

For the moment, we specify neither initial nor boundary conditions for the second equation, we just assume that corresponding problem is well posed and that it admits an $\LLd$ solution. 
For the coefficients of both the  operators  we assume that are merely bounded and continuous.  

\begin{proposition}
\label{av-obs-we}
Suppose that there is a constant $\tilde C$, time $T$ and an open subdomain $\omega$ such that for any choice of initial conditions $ \beta_0, \beta_1$  the solution $u_1$ of \eqref{discrete-wave} satisfies
\begin{equation}
\label{assumption}
E_1(0):=\nor {\beta_0}{\LLd}^2 + \nor {\beta_1}{\HH{-1}}^2
\leq \tilde C \int_0^T \int_\omega |  u_1 |^2 d\mx dt \,.
\end{equation}
In addition, we assume that characteristic sets $\{p_i(t, \mx, \tau, \mxi)=0\}, i=1,2$ have no intersection for $(t, \mx)\in\oi 0T\times\omega, (\tau, \mxi) \in \Sd$,  where $p_i$ stands for the principal symbol of the operator $P_i$.

Then for any $\theta_1, \theta_2 \in\R, \theta_1\not=0$ there exists a constant $\tilde C_{\theta_1}$ such that the observability inequality
\begin{equation}
\label{obs-ineq}
E_1(0)
\leq  \tilde C_{\theta_1} \left(\int_0^T \int_\omega |\theta_1  u_1 + \theta_2 u_2|^2 d\mx dt + \nor {\beta_0}{\HH{-1}}^2 + \nor {\beta_1}{\HH{-2}}^2\right)
\end{equation}
holds for any pair of solutions $(u_1, u_2)$ to \eqref{discrete-wave}.
\end{proposition}

{\bf Proof:} 
Of course, the case $\theta_2=0$ holds trivially, and for simplicity is excluded from the further analysis.

We argue by contradiction.
Assuming the contrary, there exists a sequence of solutions $u_1^n, u_2^n$ such that 
\begin{equation}
\label{avcontradict}
E_1^n(0)>n \left(\int_0^T \int_\omega |\theta_1  u_1^n + \theta_2  u_2^n|^2 d\mx dt
+ \nor {\beta_0^{n}}{\HH{-1}}^2 + \nor {\beta_1^{n}}{\HH{-2}}^2\right).
\end{equation}
As the considered problem is linear, without loosing generality we can assume that $E^n(0)=1$.
Thus \eqref{avcontradict} implies that   $ \nor {\beta_0^{n}}{\HH{-1}}^2 + \nor {\beta_1^{n}}{\HH{-2}}^2 \to 0$, resulting in the weak convergence   $(\beta_0^{n}, \beta_1^{n}) \dscon (0, 0)$ in $\Ld\Omega \times \H{-1}\Omega$. Therefore the  solutions  $(u_1^n)$ converge weakly to zero in $\Ld{\Omega \times \oi 0T}$ as well. In order to obtain a contradiction, we have to show that  the last convergence is strong, at least on the observability region.

From the contradictory assumption \eqref{avcontradict} we have that   the H-measure   $\nu$ associated to a  subsequence of $(\theta_1  u_1^n + \theta_2  u_2^n)$ vanishes on $\oi 0T\times\omega$. Furthermore, it is of the form
$$
\nu=\theta_1^2 \mu_1+\theta_2^2 \mu_2+ \theta_1 \theta_22\Re \mu_{12},
%\leqno(\for)
$$
where on the right hand side  the elements of the matrix measure associated to the vector subsequence of $( u_1^n,  u_2^n)$ are listed, with $\mu_{12}$ denoting the off-diagonal element. Note that $(u_2^n)$ is bounded in $\Ld{\oi 0T \times\omega}$, since that is the case for $(u_1^n)$ (by boundedness of initial data), and for the linear combination $(\theta_1  u_1^n + \theta_2  u_2^n)$ (by contradictory assumption \eqref{avcontradict}), which enables one to associate an H-measure to it.

According to the localisation property  for H-measures, each $\mu_j$ is supported within the corresponding characteristic set $\{p_i(t, \mx, \tau, \mxi)=0\}$, $i=1, 2$, which, by assumption, are disjoint on the observability region. On the other hand, from the very definition of matrix H-measures  it follows that off-diagonal entries are dominated by the corresponding diagonal elements. More precisely, it holds that 
$\supp \mu_{12} \subseteq \supp \mu_1 \cap  \supp \mu_2$,  implying that $\mu_{12}=0$ on the observability region. 

Thus we get that
$$
\nu=\theta_1^2 \mu_1+\theta_2^2\mu_2=0 \qquad {\rm on} \quad \oi 0T\times\omega\,.
$$

As  $\mu_1$ and $\mu_2$ are positive measures and $\theta_1>0$, it follows that 
$\mu_1$ vanishes on $\oi 0T\times\omega$ as well.
Thus we get strong convergence of $(u_1^n)$ in $\Ld{\oi 0T \times\omega}$, which together with the assumption of the  constant, non-zero initial energy contradicts  the observability estimate \eqref{assumption}.
\hfill $\Box$

\begin{remark}
The last result provides  surprising stability of the observability estimate \eqref{assumption} under uncertain perturbations, up to compact reminders. Essentially, the only requirement for the perturbation is separation of the characteristic sets. This implies that the wave component can be observed robustly when adding unknown perturbations, up to a finite number of low frequencies.
\end{remark}

In the next step we would like to obtain the strong observability inequality for initial energy $E_1(0)$ by removing compact terms in \eqref{obs-ineq}. To this effect, we have to specify some additional constraints on  the problem for the perturbation $u_2$.

We take $P_2$ to be an evolution operator of the form
\begin{equation}
\label{P2}
P_2=(\partial_t)^k + c_2(\mx) \lA_1, \quad k\in \N,
\end{equation}
where $\lA_1$ is an elliptic part of the wave operator $P_1$, while $c_2$ is a  bounded and continuous function.

\begin{theorem}
\label{av-obs}
As above, we assume that the coefficients of the operator $P_1$ are  bounded and continuous, and that  the corresponding solution $u_1$ satisfies the observability inequality  \eqref{assumption}.  

In addition we assume that  the perturbation operator $P_2$  is of the form \eqref{P2}. In the case $k=2$ (i.e. $P_2$ being a wave operator) the separation of coefficients $c_2(\mx)-1\not= 0$ is supposed on $\omega$.

For 
the initial values of solutions $u_i, i=1,2$ we supposed these are related by a linear operator such that whenever  $\left(( \theta_1 u_1(0)+ \theta_2 u_2(0))|_{\omega}=0\right) $ then $\left(u_1(0)|_{\omega}=u_2(0)|_{|\omega}=0\right)$, and the analogous implication holds for the initial first order time derivatives.

Then there is a positive constant $C_{\theta_1}$ such that the strong observability inequality holds:
\begin{equation}
\label{obs-ineq-str}
E_1(0)
\leq  C_{\theta_1}  \int_0^T \int_\omega |\theta_1  u_1 + \theta_2  u_2|^2 d\mx dt \,.
\end{equation}
\end{theorem}
\begin{remark}
 Note that the above assumptions directly imply that characteristic sets of  $P_1$ and $P_2$ are disjoint. Indeed, for $P_2$  being an evolution operator of order $k\not=2$, its principal symbol equals 0 only in poles $\mxi=\vnul$ (case $k=1$), or on the equator $\tau=0$ (case $k>2$)  of the unit sphere in the dual space, where $p_1=\tau^2-\mA_1(\mx)\mxi^2$ differs from zero.

In the case $k=2$ separation of the  characteristic sets is provided by the assumption $c_2(\mx)\not=1$ on $\omega$.
\end{remark}

{\bf Proof:} 
As in the proof of Proposition \ref{av-obs-we}, let us suppose the contrary. Then there exists  a sequence of solutions $u_1^n, u_2^n$ to \eqref{discrete-wave} such that $E_1^n(0)=1$ and
\begin{equation}
\label{contra2}
\int_0^T \int_\omega |\theta_1  u_1^n + \theta_2  u_2^n|^2 d\mx dt \str 0\,.
\end{equation}
Thus the corresponding weak limits satisfy both the equation $P_i u_i=0$, as well as the relation 
\begin{equation}
\label{relate}
\theta_1  u_1 + \theta_2  u_2=0 
\end{equation}
on the observability region. As functions $u_i, \partial_t u_i$ are continuous with respect to time, the bounds on initial data imply $\beta_0=\beta_1=0$ on $\omega$.

Assumptions on the operators $P_i, i=1,2$ ensure that corresponding characteristic sets do not intersect. By applying localisation property of H-measures as in the proof of Proposition \ref{av-obs-we} we get that $u_1^n$ converges to $u_1$ strongly on $\oi 0T\times\omega$.

It remains to show that the limit $u_1$ vanishes on the observability region, which, together with the assumption of the constant non-zero initial energy, will contradict the observability assumption \eqref{assumption}.

We split the rest of the proof into several cases.
\begin{itemize}
\item[a)] {$\bf (k=2)$} Due to the relation \eqref{relate} it follows that 
$$
(c_2-1) \partial_{tt} u_1=0, \qquad (t, \mx) \in \oi 0T \times \omega\,.
$$
As $|c_2-1|>0$, and the initial data are $0$ on $\omega$, it implies $u_1=0$ on the observability region.
\item[b)] {$\bf (k=1)$}
Relation \eqref{relate} implies 
$$
\partial_{t}u_1={ c_2(\mx)} \partial_{tt}u_1,
$$
which together with $u_1(0, \cdot)=\partial_{t} u_1(0, \cdot)=0$ on $\omega$ provides the claim.
\item[c)] {$\bf (k>2)$}
Similar as above we obtain 
$$
\partial_{t}^{k-2} (\partial_{tt}u_1) = c_2(\mx) \partial_{tt}u_1\,.
$$
As $\partial_{tt}u_1(0, \cdot)=- \lA_1u_1(0, \cdot)=0$ on $\omega$, and similarly for the higher order derivatives, the claim follows. 
\end{itemize}

\hfill $\Box$

\begin{remark}
\label{remark}
Several remarks are in order.
\begin{itemize}
\item The observability assumption \eqref{assumption} on a solution of the wave equation is equivalent to the Geometric Control Condition (GCC, \cite{BLR}), stating that projection of each bicharacteristic ray on a physical space has to  enter the observability region in a finite time.

\item The last theorem also holds if, instead of initial data of two components being linked by an operator,  we assume  a cone condition $\nor{u_2(0)}{\Ld\omega} \leq c \nor{u_1(0)}{\Ld\omega}$ $\left(\nor{u_2(0)}{\Ld\omega} \geq c \nor{u_1(0)}{\Ld\omega}\right)$, with a constant $c<\theta_1/\theta_2 \left(c>\theta_1/\theta_2 \right)$. The latter condition is stable under passing to a limit, and also ensures the implication  $\left(( \theta_1 u_1(0)+ \theta_2 u_2(0))|_{\omega}=0\right) \Longrightarrow \left(u_1(0)|_{\omega}=u_2(0)|_{|\omega}=0\right)$, which suffices for the proof.

%\item The assumption on initial data being related by a linear operator can be replaced by a cone condition
%$$\nor{u_2(0)}{\Ld\omega}\leq c \nor{u_1(0)}{\Ld\omega}, \quad c< \theta_1/\theta_2\,,$$
%or any other condition, stable under passage to a limit, such that the implication  $\left(( \theta_1 u_1(0)+ \theta_2 u_2(0))|_{\omega}=0\right) \\ \Longrightarrow \left(u_1(0)|_{\omega}=u_2(0)|_{|\omega}=0\right)$ holds.

\item The result \eqref{obs-ineq-str} can be generalise to a more general perturbation operator $P_2$  by assuming that coefficients of both the operators are analytic. In that case the separation of characteristic sets implies the separation of corresponding analytic wave front sets (\cite[Theorem 9.5.1]{Hor2}). Together with \eqref{relate} it provides that $u_1$ is analytic on the observability region. Constraints on initial data and finite velocity of propagation imply $u_1=0$ on an open set near $t=0$, and as the solution is analytic it vanishes on the whole  observability region which contradicts  the observability assumption \eqref{assumption}.

\item If $P_2$ is a wave operator the strong observability inequality  \eqref{obs-ineq-str} is equivalent to the controllability of a suitable linear combination (determined by the operator linking the initial data) of  solutions to the adjoint system under a single control (cf.~\cite{LZ}). 

Meanwhile, the weak observability result \eqref{obs-ineq} in that case corresponds to the average controllability of the adjoint system up to a finite number of low frequencies. 

\item The last theorem generalises the results of \cite{LZ} by allowing for a  general evolution  operator $P_2$ which does not have to be the wave one. 

In addition, it allows for an arbitrary linear combination of system components, while in \cite{LZ} just their (weighted) average is explored. Specially, if the difference $u_1-u_2$ is considered, the result corresponds to the synchronisation problem (e.g.~\cite{LRH}) in which all the components are driven to the same state by applying the null controllability of their differences.  

Furthermore, unlike in \cite{LZ}, the proof of the relaxed observability inequality \eqref{obs-ineq} does not rely on the propagation property of H-measures, which allows for system's coefficients to be merely continuous. On the other hand, such approach avoid technical issues related to the reflection of H-measures on the domain boundary.

\item The theorem to some extent also generalises the results of \cite{Z} in which a similar result is provided  for the system \eqref{discrete-wave} consisting of a wave and a heat operator with constant  coefficients (or more generally with a common elliptic part).

However, although allowing for a more general perturbation operator, it requires initial data of two components of the system \eqref{discrete-wave} to be related, while in \cite{Z} no assumptions on initial data for the second component is assumed.

\item The constant $C_{\theta_1}$ can be taken to be uniform for $\theta_1 \ge \theta_*$ and $\theta_* >0$.

\item The weak observability result \eqref{obs-ineq} is easily generalised to a system with a finite number of components, under assumption that the characteristic set of the leading operator $P_1$ is separated from the characteristic sets of all the other operators, while the latter ones can be  arbitrary related. 

However, the generalisation of the result to an infinite dimensional setting is not straightforward. It requires study of the localisation property for H-measures determined by a sequence of function series,  and is the subject of the next subsection. 

On the other side, the generalisation of  the strong estimate result \eqref{obs-ineq-str} to a system consisting of more than two components has still not been obtained, and is a subject of the current investigations. 

\item The result \eqref{obs-ineq} also holds if the observability region is not of a cylindrical type, but a more general set satisfying the GCC for the first component of the system. Such generalisation corresponds to a moving control (cf.~\cite{LRTT}). However, possible derivation of the corresponding stronger observability result remains open, due to the last part of the proof of  Theorem \ref{av-obs}, in which the special (cylindrical) shape of the observability region is used.

\item Note that the observability result  \eqref{obs-ineq-str} is weaker than the one required in the simultaneous control, where one has to estimate initial energy of all components entering the system. Of course, the assumptions in the latter case are stronger, as one has  to assume that the observability set satisfies the GCC for the second component as well.
\end{itemize}
\end{remark}

\subsection{Infinite discrete setting}
In this subsection we want to analyse the stability of the observability estimates for the wave equation when the perturbation is given as a superposition of infinitely many components, each determined by a differential operator $P_i$, which, in general, does not have to be the wave one.   Thus the system of interest reads as:
\begin{equation}
\label{infinite-discrete}
\begin{split}
P_1 u_1= \partial_{tt} u_1 - \dv (\mA_1(t, \mx) \nabla u_1)&=0, \qquad (t, \mx) \in \Rpl\times\Omega\cr
P_i u_i&=0, \qquad (t, \mx) \in \Rpl\times\Omega, i\geq 2\cr
u_1&=0, \qquad (t, \mx) \in \Rpl\times\partial\Omega\cr
u_1(0, \cdot)&=\beta_0 \in \Ld\Omega\cr
  \partial_t u_1(0, \cdot)&=\beta_1 \in \H{-1}\Omega \,,\cr
\end{split}
\end{equation}
with the same assumptions on the domain $\Omega$ and the operator $P_1$ being assumed for   the system \eqref{discrete-wave}.  For the other equations, neither initial nor boundary conditions are specified. For the moment, we just assume the corresponding coefficients are bounded and continuous, and the problems are well defined with solutions in $\Ld\Omega$. 

In this setting, the same microlocal analysis tool as in the finite case, in particular the localisation property of H-measures, is applied in the study of the stability of the observability estimates. However, as perturbations are determined by a superposition of infinitely many solutions, this requires analysis of the mentioned property  for a sequence  of function series.

Namely, it is well known that an H-measure associated to a linear combination of two sequences is supported within the union of supports of measures determined by each component, and the same property holds for any finite linear combination. However, in general it fails when considering superposition of  infinite many sequences, as shown by the next example.
\begin{example}
Let $(u^n)$ and $(f^n)$ be $\Ld\Rd$ sequences, whose corresponding H-measures $\mu_u$ and $\mu_f$ have disjoint supports, and let $(\theta_i)$ be a sequence of nonnegative numbers summing in 1.

Define the following sequences
$$
v_i^n=\left\{
\begin{array}{cc}
\theta_i u^n & i \not= n\\
f^i & i=n\,.\
\end{array}
\right.
$$
Thus for each $i$ an H-measure $\nu_i$ associated to $v_i^n$ equals $\theta_i^2 \mu_u$.

On the other side we have that $\sum_i v_i^n= (1-\theta_n) u^n + f^n$, and the corresponding measure equals $\mu_u+\mu_f$. 

\end{example}

Thus in order to constrain support of an H-measure by supports of corresponding components we have to impose additional assumptions on constituting sequences. More precisely, the following result holds.
\begin{lemma}
\label{local-lemma}
Let $(\theta_i)$ be an averaging sequence of positive numbers summing to 1, and let $(u_i^n)_n, i\in \N$ be a family of uniformly bounded $\LLd$ sequences, i.e. we assume there exists a constant $C_u$ such that $\nor{u_i^n}\LLd \leq C_u, i,n \in \N$. 

Define the linear combination $v_n= \sum_i \theta_i u_i^n$, and denote by 
 $\mu_i$ and $\nu$   H-measures associated to (sub)se\-quen\-ces (of) $(u_i^n)_n$ and $(v_n)$, respectively. Then 
\begin{equation}
\label{local-result}
\supp \nu \subseteq {\rm Cl}\Big(\cup_{i} \supp \mu_i\Big).
\end{equation}
\end{lemma}
{\bf Proof:} 
Take an arbitrary pseudodifferential operator of order zero, $P\in \Psi_c^0$, with a symbol $p( \mx, \mxi)$ being compactly supported within the complement of the closure of $\cup_{i} \supp \mu_i$. 

By the definition of H-measures we have
\begin{equation}
\label{H-series}
\Dup {\nu}p
=\lim_n \int_{\Rd} P \Big( \sum_1^\infty \theta_i  u_i^n\Big)( \mx) 
\Big(\sum_1^\infty \theta_j  u_j^n\Big)( \mx)  d\mx\,.
\end{equation}

As $P$ is a continuous operator on $\Ld\Rd$ it follows that 
\begin{equation}
\begin{split}
\lim_n \!\left| \int_{\Rd}   P \Big(\! \sum_k^\infty \theta_i  u_i^n\Big)( \mx) 
\Big(\!\sum_1^\infty \theta_j  u_j^n\Big) ( \mx)  d\mx\right|
& \leq \limsup_n C_P \Bnor{\sum_k^\infty \theta_i  u_i^n}{\LLd} 
\Bnor{\sum_1^\infty \theta_i  u_i^n}{\LLd} \cr
&\leq \limsup_n C_P \,C_u^2 \left(\sum_k^\infty \theta_i \right) \buildrel k \over \dstr 0 \,, \cr
\end{split}
\end{equation}
where $C_P$ is the $\LLd$ bound of the operator $P$. The last sum is a remainder of a convergent series, and the above limit  converges to zero uniformly with respect to $n$.

Similarly, one shows the same property holds for  $\lim_n \left| \int_{\Rd}   P \Big( \sum_1^\infty \theta_i^n  u_i^n\Big)( \mx) 
\Big(\sum_l^\infty \theta_j^n  u_j^n\Big) ( \mx)  d\mx\right|$. 
Thus we can exchange limits in \eqref{H-series}, getting
\begin{equation}
\begin{split}
\Dup {\nu}p
&=\sum_{i=1}^\infty\sum_{j=1}^\infty \lim_n \int_{\Rd} P \left(  \theta_i  u_i^n\right)(\mx)\, \theta_j  u_j^n  ( \mx)  d\mx\cr
&=\sum_{i=1}^\infty\sum_{j=1}^\infty \theta_i\theta_j \Dup {\mu_{ij}}p=0,
\end{split}
\end{equation}
where $\mu_{ij}$ are H-measures determined by sequences $(u_i^n)$ and $(u_j^n)$, supported within the closure given in \eqref{local-result}, outside which $p$ is supported.
\hfill $\Box$

As a consequence of the last lemma, in order to  apply the localisation property  within the analysis of  observability estimates for solutions to \eqref{infinite-discrete},  some kind of uniform boundedness on the solutions has to be assumed.

\begin{proposition}
\label{weak-infinite}

Suppose the observability inequality \eqref{assumption} holds for a solution $u_1$ to \eqref{infinite-discrete}.  

As for the system \eqref{infinite-discrete}, suppose that $\LLd$ norm of all the  solutions $u_i$  is  dominated (up to a multiplicative constant, independent of a choice of initial data) by the energy norm of $u_1$.
In addition assume that characteristic set $\{p_1(t, \mx, \tau, \mxi)=0\}$ has no intersection with ${\rm Cl}\Big(\cup_{i\geq2} \{p_i(t, \mx, \tau, \mxi)=0\}\Big)$ for $(t, \mx)\in\oi 0T\times\omega, (\tau, \mxi) \in \Sd$,  where $p_i$ stands for the principal symbol of the operator $P_i$.

Then for any  averaging sequence $(\theta_i)$ of positive numbers summing to 1, with $\theta_1>0$, there exists a constant $\tilde C_\theta$ such that the observability inequality
\begin{equation}
\label{obs-ineq-infinite}
E_1(0)
\leq  \tilde C_\theta \left(\int_0^T \int_\omega |\sum \theta_i u_i|^2 d\mx dt + \nor {\beta_0}{\HH{-1}}^2 + \nor {\beta_1}{\HH{-2}}^2\right)
\end{equation}
holds for any family of solutions $(u_i)$ to \eqref{infinite-discrete}.
\end{proposition}

{\bf Proof:} 
Assume the contrary.
Then there exist sequences of initial conditions $(\beta_0^n), (\beta_1^n)$, and of associated solutions $(u_i^n)$, such that  
\begin{equation}
\label{avcontradiction}
1=E_1^n(0):=\nor{ \beta_0^n}\LLd^2 + \nor{\beta_1^n}{\HH{-1}}^2>n \left(\int_0^T \int_\omega |\sum_{i=1}^\infty \theta_i  u_i^n |^2 d\mx dt + \nor { \beta_0^n}{\HH{-1}}^2 + \nor {\beta_1^n}{\HH{-2}}^2\right).
\end{equation}
Let $\nu$ be  an H-measure associated to a (sub)sequence of $\sum_{i=1}^\infty \theta_i^n u_i^n$. 
Due to the inequality \eqref{avcontradiction}, it equals zero on $\oi 0T\times\omega$. 

We split the last sum into two parts $\theta_1  u_1^n + \sum_{i=2}^\infty \theta_i  u_i^n$, and we  rewrite $\nu$ in the form
$$
\nu= \nu_1+ \nu_2+ 2\Re \nu_{12},
%\leqno(\for)
$$
where $ \nu_1$ and $\nu_2$ are H-measures associated to  (sub)sequences (of) 
$(\theta_1^n  u_1^n)$ and $(\sum_2^\infty \theta_i^n  u_i^n)$, respectively, while 
$\nu_{12}$ is a measure corresponding to their product. In addition, $\nu_1=\theta_1^2 \mu_1$, where by $\mu_i$ we denote a  measure  associated to a  (sub)sequence (of) the $i-$th component $u_i^n$. 

From here the statement of the theorem is obtained easily  (following the lines of the proof in finite discrete case, Proposition \ref{av-obs-we}), once we show that $\nu_1$ and $\nu_2$ have disjoint supports.

By the localisation property for H-measures, each measure $\mu_i$ is supported within the set $\{p_i(t, \mx, \tau, \mxi)=0\}$.

The assumption on the domination of solutions to \eqref{infinite-discrete} by an energy norm of $u_1$, together with the constant initial energy  $E_1^n(0)$ implies uniform bound on solutions $u_i^n$, both with respect to $i$ and $n$. Thus we can apply Lemma \ref{local-lemma} to conclude that $\nu_2$ is supported within the set 
$$
\label{avsupport}
{\rm Cl}\Big(\cup_{i\geq 2} \{p_i(t, \mx, \tau, \mxi)=0\}\Big),
$$
which, due to the assumption  on separation of the characteristics set, does not intersect the support of $\nu_1=\theta_1^2 \mu_1$. 
As $\theta_1$ is strictly positive, we get that $u_1^n$ converges to 0 strongly in $\Ld{\oi 0T\times \omega}$, which contradicts the observability estimate \eqref{assumption}.
\hfill $\Box$

\begin{remark} 
\label{remark2} 
$ \newline$
\begin{itemize}
\item
\vskip -5mm
The assumption of the last proposition requiring solutions $u_i$ of \eqref{infinite-discrete} to be dominated by the energy norm of $u_1$  occurs, for example, in a case of  a system consisting of the operators  of the same form, $P_i=\tau^2 - \mA_i(t, \mx) \mxi \cdot \mxi$, with uniformly bounded (both from below and above) coefficients and initial energies.

The assumption on separation of characteristics sets in that case can be stated as 
$$
 \mA_1(t, \mx)\mxi \cdot \mxi \; >(<)\; \sup_{i\geq 2} \,(\inf_{i\geq 2})\; \mA_i(t, \mx)\mxi \cdot \mxi, \qquad (t, \mx) \in \oi 0T\times\omega, \; \mxi  \not= \vnul \,,
$$ 
i.e. the fastest (or the slowest) velocity is strictly separated from all the others.

In that case the weak observability \eqref{obs-ineq-infinite} result is equivalent to the averaged controllability of the adjoint system up to a finite number of low frequencies.

Of course, one can construct more general systems, including operators of different types as well, that satisfy the required boundedness assumption.
\item As already mentioned in previous subsection,  obtaining corresponding strong observability result in this setting remains an  open problem. 
\item
The constant $C_\theta$ appearing in \eqref{obs-ineq-infinite} can be taken uniformly for a family of averaging sequences, each   satisfying
\begin{itemize}
\item[({\it i})] $\theta_1\geq \theta_\ast$,
\item[({\it ii})] $\sum_k^\infty \theta_i \leq \epsilon_k$,
\end{itemize}
where $\theta_\ast\in\ozi 01$  and  $(\eps_k)$ is  a null sequence, both independent of a choice of a particular sequence  $(\theta_i)$.
\end{itemize}
 \end{remark}

\subsection{Simultaneous observability}
\label{Simult-subsection}
The subsection deals with a problem of recovering energy of a system by observing an average of solutions on a suitable subdomain. For this purpose one has to estimate initial energies of all system components, unlike the case of the average observability where this was required just for the first one. A two component system is analysed firstly, while  generalisations to a more dimensional case is discussed at the end. 

We reconsider the system \eqref{discrete-wave}  assuming that $P_2$ is an evolution operator of the form
\begin{equation}
\label{P2-s}
P_2=(\partial_t)^k + \lA_2, \quad k\in \N,
\end{equation}
where $\lA_2$ is an (uniformly) elliptic operator (in general different from $\lA_1)$, and the problem for  the perturbation $u_2$ is accompanied by a series of initial conditions
$$
\left((\partial_t)^j u_2\right)(0)= \gamma_j \in \H{-j}{\Omega},\quad j=0, \dots, k-1.
$$
Its initial energy is denoted by 
$$
E_2(0)=\sum_0^{k-1} \nor{\gamma_j}{\H{-j}{\Omega}}\,.
$$

As in the previous subsection, we start with a weak observability result.
\begin{proposition}
\label{sim-obs-we}
Suppose that there is a constant $\tilde C$, time $T$ and an open subdomain $\omega$ such that for any choice of initial conditions  the solutions  to \eqref{discrete-wave} satisfy
\begin{equation}
\label{assumption-sim}
E_i(0)\leq \tilde C \int_0^T \int_\omega |  u_i |^2 d\mx dt, \quad i=1,2 \,.
\end{equation}
In addition assume that characteristic sets $\{p_i(t, \mx, \tau, \mxi)=0\}, i=1,2$ have no intersection for $(t, \mx)\in\oi 0T\times\omega, (\tau, \mxi) \in \Sd$,  where $p_i$ stands for the principal symbol of the operator $P_i$.

Then for any $\theta_1, \theta_2 \in\R\setminus\{0\}$ there exists a constant $\tilde C_{\theta}$ such that the observability inequality
\begin{equation}
\label{obs-ineq-sim}
E_1(0)+E_2(0)
\leq  \tilde C_{\theta} \left(\!\int_0^T \!\!\int_\omega |\theta_1  u_1 + \theta_2 u_2|^2 d\mx dt + \nor {\beta_0}{\HH{-1}}^2 + \nor {\beta_1}{\HH{-2}}^2 
+ \nor {\gamma_0}{\HH{-1}}^2 +\ldots+ \nor {\gamma_{k-1}}{\HH{-k}}^2\!\right)
\end{equation}
holds for any pair of solutions $(u_1, u_2)$ to \eqref{discrete-wave}.
\end{proposition}
The result is obtained easily by following the steps of the proof presented above in the averaged observability setting. Assuming the contrary and implying microlocal analysis tools, one shows that both components $u_i^n$ converge to 0 strongly on the observability region, thus obtaining the contradiction. 

However, a different approach is required in  order to obtain the strong observability inequality for initial energy  by removing compact terms in \eqref{obs-ineq-sim}. It is based on a standard compactness-uniqueness procedure of reducing the observability for low frequencies to an elliptic unique continuation result \cite{BLR, DLR}.

We introduce  a subspace $N(T)$ of $H=\Ld{\Omega}\times \H{-1}{\Omega}\times\Ld{\Omega}\times \dots \times \H{1-k}{\Omega}$, consisting of initial data for which the average of solutions to \eqref{discrete-wave} vanishes on the observability region
\begin{equation*}
N(T):=\{ (\beta_0, \beta_1, \gamma_0, \dots, \gamma_{1-k}) \in H |\,  \theta_1 u_1+ \theta_2 u_2 = 0 \; {\rm on}\; \oi 0T \times \omega\}.
\end{equation*}
Based on the relaxed observability inequality \eqref{obs-ineq-sim}  it follows that $N(T)$ is a finite dimensional space. Furthermore, the following characterisation holds.

\begin{lemma}
\label{uc-lemma}
We assume one of the following statements holds:
\begin{itemize}
\item[a)] The order $k$ of time derivative in \eqref{P2-s}  is odd. Coefficients of both the operators $P_1$ and $P_2$ are
time independent and of class $\CC^{1,1}$,
\item[b)] The time derivative order $k$ is even, and $ \lA_1^{k/2} - \lA_2$ (or $ -(\lA_1^{k/2} - \lA_2)$) is an uniformly elliptic operator. 
Coefficients of both the operators $P_1$ and $P_2$ are analytic. 
\end{itemize}
Then 
$
N(T) = \{\vnul \}.
$
\end{lemma}
{\bf Proof:} 
One first shows that $N(T) $ is an $A$-invariant,  where $A$ is an unbounded operator on $H$:
\begin{equation}
\label{A}
A= \left(\begin{array}{c|c}
\begin{matrix}
\mnul& -1&\\
\lA_1& \; \mnul\\
\end{matrix}
& \mnul\\
\hline
\mnul&
\begin{matrix}
\mnul&
-\mI_{k-1}\\
\lA_2&\; \mnul\\
\end{matrix}\\
\end{array}
\right)\,,
\end{equation}
with the domain 
$D(A)=  \Hnl{1}{\Omega}\times \Ld{\Omega}\times\Hnl{1}{\Omega} \times \cdots \times \H{2-k}{\Omega}$.

Being $A$-invariant and finite-dimensional, it contains an eigenfunction of $A$. Thus there is a $\lambda\in \C$ and $ (\beta_0, \beta_1, \gamma_0, \dots, \gamma_{1-k}) \in N(T)$ such that 
\begin{equation}
\label{eigenpair}
\begin{split}
\lA_1 \beta_0&=- \lambda^2 \beta_0\cr
\lA_2 \gamma_0&=(-1)^{k-1}\lambda^k  \gamma_0\cr
\beta_1&=-\lambda \beta_0\cr
\gamma_j&=(-1)^{j}\lambda^j  \gamma_0,\quad  j=1, \dots, k-1\,.\cr
\end{split}
\end{equation}
By the definition of $N(T) $ it follows $\theta_1 u_1+ \theta_2 u_2=0$ on $\oi 0T \times \omega$, and specially 
\begin{equation}
\label{vanish}
\theta_1 \beta_0+\theta_2 \gamma_0=0\quad {\rm on} \quad\omega\,.
\end{equation}
At this level, we want to show that each  assumption of the lemma implies $\beta_0=\gamma_0=0$.
\begin{itemize}
\item[a)] 
As $\lA_i, i=1,2$ are positive operators and $k$ is odd, from \eqref{eigenpair} it follows that one of functions  $\beta_0, \gamma_0$ is trivial. By relation \eqref{vanish} it follows that the other one also equals zero on $\omega$. Being an eigenfunction of an elliptic operator,   the unique continuation argument  (e.g.~\cite[Theorem 3]{Hor}) implies it is  zero  everywhere. 

\item[b)] Analyticity of coefficients implies analyticity of eigenfunctions. Specially it follows $\theta_1 \beta_0+\theta_2 \gamma_0=0$ everywhere, and relations \eqref{eigenpair} imply 
$$
(\lA_1^{k/2} - \lA_2)\beta_0  =0.
$$
 Assumptions on the operator $\lA_1^{k/2} - \lA_2 $  imply $\beta_0=0$ on $\Omega$. 
\end{itemize}
\hfill $\Box$
\begin{remark}
\label{remark-lemma}
In a special case $\lA_2=-\dv (c_2(\mx) \nabla)$ one easily proves that the last Lemma holds with analytic coefficients $c_1, c_2$ being separated just on an arbitrary non-empty open set, and not on the entire $\Omega$.  
\end{remark}

\begin{theorem}
Under the assumptions of Proposition \ref{av-obs-we} and Lemma \ref{uc-lemma} there is a positive constant $C_\theta$ such that the strong observability inequality holds:
\begin{equation}
\label{obs-ineq-strong-sim}
E_1(0)+E_2(0)
\leq  C_\theta  \int_0^T \int_\omega |\theta_1  u_1 + \theta_2  u_2|^2 d\mx dt \,.
\end{equation}
\end{theorem}
{\bf Proof:} 
As in the proof of Proposition \ref{av-obs-we}, let us suppose the contrary. Then there exists  a sequence of solutions $u_1^n, u_2^n$ to \eqref{discrete-wave} such that $E_1^n(0)=1$ and 
$$
\int_0^T \int_\omega |\theta_1  u_1^n + \theta_2  u_2^n|^2 d\mx dt \str 0\,.
$$
Thus for weak limits $(u_1, u_2)$ of solutions on the observability region we have $\theta_1  u_1 + \theta_2  u_2=0$, implying $\left(u_1(0), \partial_t u_1(0), u_2(0), \dots, \left((\partial_t)^{k-1} u_2\right)(0)\right)\in N(T)$. By means of the above lemma and taking into account the relaxed observability inequality, it follows
$$
1\leq \tilde C_\theta \left(\int_0^T \int_\omega |\theta_1  u_1^n  + \theta_2  u_2^n |^2 d\mx dt 
+ \nor {\beta_0^{n}}{\HH{-1}}^2 + \nor {\beta_1^{n}}{\HH{-2}}^2
+ \nor {\gamma_0^{n}}{\HH{-1}}^2 +\ldots+ \nor {\gamma_{k-1}^{n}}{\HH{-k}}^2\right) \str 0,
$$
thus obtaining a contradiction.
\hfill $\Box$
\vskip 2mm

We close this subsection by the following remarks.
\begin{remark}
\label{remark-sim}
$\newline$
\begin{itemize}
\item 
\vskip -5mm

If $P_2$ is a second order evolution operator the strong observability inequality  \eqref{obs-ineq-strong-sim} is equivalent to the simultaneous  controllability  of the adjoint system, also studied in \cite{LZ}, by which one controls each component individually (and not just their average). 
\item The notion of simultaneous observability is stronger than the average one, as it estimates energy of all system components, whose initial data, in this case, are not related. Consequently, it requires stronger assumption of GCC being satisfied by each component. 
\item The application of the compactness-uniqueness procedure in the passage from the weak to the strong observability estimate allows the perturbation $P_2$ to be an evolution operator with an arbitrary elliptic part. However, such approach is not possible in the averaged observability setting. Namely, in order for subspace $N(T)$ to be finite dimensional one has to relate the initial data of two components by a bounded linear operator. But such  constrain would not be preserved under action of the operator $A$ given by \eqref{A}, and as a consequence $N(T)$ would not be $A$-invariant. 
\item The weak observability result \eqref{obs-ineq-sim} is easily generalised to a system with a finite number of components, under assumption that the characteristic sets of all  operators are mutually disjoint.

\item As in the averaged observability case, the  generalisation of strong estimate result \eqref{obs-ineq-strong-sim} to a system consisting of more than two components has still not been obtained, and is a subject of the current investigations.

\end{itemize}
\end{remark}

%\subsection{relation to the control theory}

\section{Observation of the Schr\"odinger equation under uncertain perturbations}
\label{section-Schrod}

In this section we consider  a system in which the first component, the one whose energy is observed, satisfies the Schr\"odinger equation, while the second one, corresponding to a perturbation, is governed by an evolution  operator $P_2$: 
\begin{equation}
\label{scrod-sys}
\begin{split}
P_1 u_1= i \partial_{t} u_1 + \dv (\mA_1(\mx) \nabla u_1)&=0, \qquad (t, \mx) \in \Rpl\times\Omega\cr
P_2 u_2&=0, \qquad (t, \mx) \in \Rpl\times\Omega\cr
u_1&=0, \qquad (t, \mx) \in \Rpl\times\partial\Omega\cr
u_1(0, \cdot)&=\beta_0 \in \Ld\Omega \,.\cr
\end{split}
\end{equation}
As in the study of perturbations of the wave dynamics in Section \ref{wave}, we specify no initial or boundary conditions for the second operators, we just assume that corresponding problem is well posed and that it admits an $\LLd$ solution. As for the system coefficients, as before we impose  merely boundedness and continuous assumptions, and suppose that $\mA_1$ is a positive definite matrix field.

\subsection{Averaged observability under non-parabolic perturbations}
\label{Schrod-hyper}
For the reasons  explained below, in this subsection we restrict the analysis to evolution  operators $P_2$ of order strictly larger than one. In that case the stability of the  Schr\"odinger observability estimate is given by the next theorem.
\begin{theorem}
\label{scrod-result1}
Suppose that there is a constant $\tilde C$, time $T$ and an open subdomain $\omega$ such that for any choice of initial datum $\beta_0$  the solution $u_1$ of \eqref{discrete-wave} satisfies
\begin{equation}
\label{assumption-schrod}
E_1(0):=\nor {\beta_0}{\LLd}^2 
\leq \tilde C \int_0^T \int_\omega |  u_1 |^2 d\mx dt \,.
\end{equation}
In addition, for the system \eqref{scrod-sys} we assume  the following:
\begin{itemize}
\item[a)] The perturbation operator is  an evolution operator of the form \eqref{P2} and of the order $k>1$. 
\item[b)] The initial values of solutions $u_i, i=1,2$  are related by a linear operator such that whenever  $\left(( \theta_1 u_1(0)+ \theta_2 u_2(0))|_{\omega}=0\right) $ then $\left(u_1(0)|_{\omega}=u_2(0)|_{|\omega}=0\right)$.
\end{itemize}

Then, for any $\theta_1 \in \ozi 01$ there exists a constant $C_\theta$ such that the observability inequality
\begin{equation}
\label{obs-ineq-scrod}
E_1(0)
\leq C_\theta \int_0^T \int_\omega |\theta_1  u_1 + \theta_2 u_2|^2 d\mx dt 
\end{equation}
holds for any pair of solutions $(u_1, u_2)$ to \eqref{discrete-wave}.
\end{theorem}
The proof goes similarly as for the observations of the wave equation. 
Required conditions a), b) are necessary for obtaining the strong observability inequality, without a compact term. 

On the other hand, in order to obtain  a relaxed inequality with a compact term, no assumption is required at all. 
Namely, the  assumption on separation of characteristic sets $\{p_i(t, \mx, \tau, \mxi)=0\}$ required in Proposition \ref{av-obs-we} becomes superfluous in this setting,  as being directly satisfied by an arbitrary evolution operator $P_2$ of order $k$ strictly larger than 1. Namely, its characteristic set does not contain the poles $\mxi=\vnul$ of the unit sphere in the dual space, which constitute the characteristic set of the  Schr\"odinger operator $P_1$. 

However, every  Schr\"odinger or the heat operator $P_2$ fails to satisfy the assumption. Namely, no matter the coefficients entering the equation, both have characteristic set within the  poles $\mxi=\vnul$, same as $P_1$, and localisation principle fail to distinguish corresponding H-measures. 
To analyse such a system one needs a microlocal tool better adopted to a study of parabolic problems. Namely, original H-measures were constructed with the aim of analysing  hyperbolic problems and are not capable to distinguish differences between the time and space variables that are intrinsic to parabolic equations. Their variant, parabolic H-measures, was recently introduced with the purpose of overcoming the mentioned constraint. 

\subsection{Parabolic H-measures}

Parabolic H-measures were first introduced  in \cite{ALjma}, while more exhausted introduction can be found in \cite{ALjfa}, elaborating in particular their basic properties: localisation and propagation one. The former will be used in the next subsection for proving the stability of the observability estimate for the 
Schr\"odinger operator under parabolic perturbation.

Here we expose the basic result on parabolic H-measures used in the note. 

The main idea in their construction is to replace the projection along the straight rays in the dual space, determined by the term $\mxi/|\mxi|$  in the definition of the original H-measures,  by  the one going along meridians of paraboloids $\tau=a \mxi^2$. The hypersurface on which the dual space (except the origin) is projected is a rotational ellipsoid
$$
\Pd: \dots \tau^2+{\mxi^2 \over 2}=1\,.
$$
Although the ellipsoid might seem as  an unnatural choice of the surface on which one want to construct the parabolic H-measures, crucial in its choice was that   the curves along which the projections are taken 
intersect it in the normal direction,
as it was in the classical case, where the rays radiating from origin are perpendicular to the unit sphere. The mentioned normality property enables study of propagation properties of the measures.

Concerning the applications, the most important is that the new tool  is also a kind of defect measure, in the sense that null parabolic H-measure is equivalent to strong convergence of $(\vu_{n})$ in $\Ldl{\Rjpd}$. 

In order to formulate the localisation principle, we first introduce some special anisotropic (Sobolev) function spaces 
$$
\H{{s \over 2}, \,s}{\Rjpd} := \Bigl\{ u\in\lS': k_p^{s} \hat u \in \Ld{\Rjpd}\Bigr\} \,, s \in \R,
$$
where $
k_p(\tau,\mxi) := \root 4\of{1+(2\pi\tau)^2+(2\pi|\mxi|)^4} 
$ is the weight function. 
These are Hilbert spaces and they are particular examples of more general H\"ormander spaces $B_{p,k}$ described in \cite[\paragraf 10.1]{Hor2}. 

We also  define the fractional derivative:
$\sqrt \partial_t$ as a pseudodifferential operator with a polyhomogeneous symbol $\sqrt{2\pi i \tau}$, i.e. 
$$
\sqrt \partial_t u = \overline\lF\left(\sqrt{2\pi i  \tau} \, \hat u(\tau)\right).
$$ 

\begin{theorem}
\label{local-parabolic}
{\bf (localisation principle for parabolic H-measures, \cite{ALjfa})}
Let $(\vu_n)$ be a sequence of functions uniformly compactly supported in $t$ and converging
weakly to zero in $\Ld{\Rjpd;\C^r}$, and let for $s\in \N$
$$
\sqrt{\partial_t}^s (\mA^0\vu_n) 
+ \sum_{|\malpha|=s} \partial^\malpha_\mx(\mA^\malpha\vu_n) 
\str 0 \quad \hbox{strongly in}\quad  \Hl{-{s \over 2},-s}{\Rjpd}\;,
$$
where $\mA^0, \mA^\malpha$ are continuous and bounded matrix coefficients,  while $\malpha\in\mN^d_0$.

Then for the associated parabolic H-measure  $\mmu$ we have
$$
\biggl((\sqrt{2\pi i \tau})^s \mA_0 
+ \sum_{|\malpha|=s} (2\pi i  \mxi)^\malpha \mA_\malpha \biggr) \mmu^\top
=\mnul.
$$
\vskip-4mm
\end{theorem}
In particular, the principle implies that the measure $\mmu$ is supported within the {\it parabolic characteristic set}:
$$
\det\Big( (\sqrt{2\pi i \tau})^s \mA_0 
 + \sum_{|\malpha|=s} (2\pi i  \mxi)^\malpha \mA_\malpha\Big) =0, 
\quad (\tau, \mxi)\in \Pd
$$

\begin{example} 
\label{example}
{\bf (Application of the localisation principle to various equations)}

In all the examples we assume that the equation coefficients satisfy the assumptions of the preceding theorem, i.e are continuous and bounded.
\begin{itemize}
\item {\bf The Schr\"odinger equation}

\noindent
Let $(u_n)$ be a  sequence of solutions to the Schr\"odinger equation
$$
i \partial_{t} u_n + \dv (\mA(t,\mx) \nabla u_n)=0.
$$
If $(u_n)$ is bounded in $\Ld{\Rpl\times\Rd}$ then the associated parabolic H-measure $\mu$ satisfies 
$$
(2\pi \tau + 4\pi^2 \mA(t,\mx) \mxi \cdot \mxi) \mu =0,
$$
implying $\mu$ is supported in points of the form $2\pi \tau = - 4\pi^2 \mA(t,\mx) \mxi \cdot \mxi$. 

\item {\bf The heat equation}

Let $(u_n)$ be a  sequence of solutions to the heat equation
$$
 \partial_{t} u_n - \dv (\mA(t,\mx) \nabla u_n)=0\,
$$
where  $\mA$ is a bounded, positive definite matrix field.
If $(u_n)$ is bounded in $\Ld{\Rpl\times\Rd}$ then the associated parabolic H-measure $\mu$ satisfies 
$$
(2\pi i \tau + 4\pi^2 \mA(t,\mx) \mxi \cdot \mxi) \mu =0.
$$
As (parabolic) H-measures live on a hypersurface in the dual space excluding the origin, and $\mA$ is  positive definite,  the term in the braces above never equals zero, implying $\mu$ is a trivial (null) measure.

\item {\bf The wave equation}

\noindent
Let $(u_n)$ be a  sequence of solutions to the wave equation
$$
\partial_{tt} u_n - \dv (\mA(t,\mx) \nabla u_n)=0.
$$
If $(u_n)$ is bounded in $\Ld{\Rpl\times\Rd}$ then the associated parabolic H-measure $\mu$ satisfies 
$$
 4\pi^2 \tau ^2  \mu =0,
$$
implying $\mu$ is supported in the equator $\tau=0$ of the hypersurface $\Pd$. 

\end{itemize}
\end{example}
From the given examples it is clear that by taking two Schr\"odinger operators with separated coefficients we shall be able to distinguish corresponding parabolic H-measures. Of course, the distinction is also possible if the considered operators are of different type (e.g.  the Schr\"odinger and the wave one), as it was the case with the original H-measures. This enables the generalisation of  Theorem \ref{scrod-result1} represented in the next subsection. 

\subsection{Averaged observability  under parabolic perturbations}
We reconsider  the system \eqref{scrod-sys} with the aim of obtaining the stability of the observability estimate for the 
Schr\"odinger operator under perturbation determined by the operator $P_2$ under minimal conditions on the latter. In particular, we want allow it to be the Schr\"odinger, or  some other parabolic type operator, like $P_1$.

The results of previous subsection, in particular the localisation property (Theorem \ref{local-parabolic}) provide the following result.

\begin{theorem}
\label{obs-schrod2}
In addition to assumptions of Theorem \ref{scrod-result1}  we allow $P_2$ to be an evolution operator of order one, and assume the separation of coefficients $c_2(\mx)\not= -i$ holds on $\omega$ in the case $k=1$.

Then for any $\theta_1 \in \ozi 01$ there exists a constant $C_\theta$ such that the strong observability inequality
\begin{equation}
\label{obs-schrod}
E_1(0)
\leq  C_\theta \int_0^T \int_\omega |\theta_1  u_1 + \theta_2 u_2|^2 d\mx dt
\end{equation}
holds for any pair of solutions $(u_1, u_2)$ to \eqref{scrod-sys}.
\end{theorem}

The proof goes along the same lines as the one of the Theorem \ref{av-obs}, using that the parabolic H-measures share the basic properties (positive definiteness, diagonal domination) with their original (hyperbolic) counterparts. 
Crucial in the proof is separation of corresponding parabolic characteristic sets. For an evolution operator $P_2$ of order $k\geq2$ that set is restricted to the equator $\tau=0$ of the hypersurface $\Pd$, where $p_1= 2\pi \tau + 4\pi^2 \mA(t,\mx) \mxi \cdot \mxi$ differs from zero (cf. Example \ref{example}). For $k=1$ the separation follows by the assumption $c_2(\mx)\not= -i$. 

The novelty obtained by the application of parabolic measures is that both the operators entering the system \eqref{scrod-sys} are allowed to be of the same type (e.g. two Schr\" odinger operators with coefficients being separated on the observability region $\oi 0T \times \omega$). In that case the observability estimate \eqref{obs-schrod} is equivalent to the averaged control of the adjoint system.

Of course, operators of different types are admissible as well, thus the last theorem incorporates the results of Subsection \ref{Schrod-hyper} as well.

Finally, let us mention that as in Section \ref{wave} one can obtain analogous result for the robust observability of the Schr\" odinger equation in the simultaneous and infinite discrete setting.

\section{Conclusion}

One of the most interesting aspects of the above results refers to a system determined by two operators of a same type - either the Schr\" odinger or the wave one. The weak observability estimates in that case correspond to the average control of solutions to adjoint system up to a finite number of low frequencies. The required assumptions are optimised: the coefficients are merely continuous and separated on the observability region. 
The strong averaged observability result, corresponding to the exact averaged controllability, assumes in addition  initial data of two components to be related in an appropriate manner. 
Its proof, unlike in \cite{LZ}, employs neither propagation property of H-measures nor unique continuation procedure, therefore  not requiring additional smoothness assumptions. 

The paper restricts to decoupled systems, but one could analyse more general ones with coupling in lower order terms (e.g.~\cite{LRH}). Such terms do not effect microlocal properties of solutions, thus enabling  a generalisation of the obtained estimates.

Surely of interest would be a generalisation of the results to a para\-me\-ter dependent system with a para\-meter ranging over a  continuous set. A problem has been  analysed in \cite{LuZ}, exploring the heat and Schr\" odinger equations with a randomly distributed parameter. 
The constant coefficients operators are considered, with eigenfunctions of corresponding elliptic parts being independent of the parameter value.
%The corresponding elliptic parts with constant coefficients, whose eigenfunctions are independent of the parameter value. 
Thus the explored averages   are presented as a solution (or a superposition of two solutions) to a similar evolution problem(s), which is crucial in the proof. 

As the next step in that direction, one would consider a
system of equations determined by evolution operators whose elliptic parts coincide up to a scalar function: 
$$
P(\nu)=\partial_t^k + c(\mx, \nu) \lA\,,
$$
with $\nu$ being a parameter, and $\lA$ the elliptic part of $P$. The variable dependence of the coefficient $c$ would require another techniques from those applied in \cite{LuZ}. In this paper  we have obtained a corresponding result for a parameter ranging over an infinite discrete set, but only at the level of a weak observability estimate with compact terms. 

By using transmutation techniques developed in \cite{EZ}, the simultaneous observability result of subsection \ref{Simult-subsection} can be employed in order to obtain controllability and observability properties for a system of heat equations (cf.~\cite{LZ}). Similar transmutation procedure can be constructed, transferring Schr\"odinger into wave type problems, and vice versa. Its application will result in observability estimates for a  system of Schr\" odinger equations, derived from  the corresponding results for the wave ones. It would be interesting to compare such obtained results with those obtained directly in Section \ref{section-Schrod} by means of parabolic H-measures, and to compare the efficiency of methods applied by each approach. 

\section*{Acknowledgements}
The  paper is  supported in part by Croatian Science Foundation under the project 9780. The author acknowledges E. Zuazua for inspiring and fruitful discussions on the subject.


\begin{thebibliography}{00}
% please try to use the bibitem system -
% the references should be in alphabetical order of authors' names.
% Articles with a single author first, author will 1 co-author next,
% then author with several co-authors;


% \bibitem{label}
% Text of bibliographic item


\bibitem{ALjma}
{\sc N.\,Antoni\'c and M.\,Lazar}:
 H-measures and variants applied to parabolic equations.
{\it J. Math. Anal. Appl.} {\bf 343(1)} (2008) 207--225.


\bibitem{ALjfa}
{\sc N.\,Antoni\'c and M.\,Lazar}:
Parabolic  H-measures.
{\it J. Funct. Anal.}~{\bf 265(7)} (2013) 1190--1239.

\bibitem{BLR} {\sc C.\,Bardos, G.\,Lebeau and  J.\,Rauch},  Sharp sufficient conditions for the observation, control, and stabilization of waves from the boundary. {\it SIAM J. Control Optim}. {\bf 30(5)} (1992)  1024--1065. 

\bibitem{BG97} {\sc N.\,Burq and P.\,G\'erard},  Condition n\'ecessaire et suffisante pour la contr\^olabilit\'e exacte des
 ondes.  {\it C. R. Acad. Sci. Paris S\'er. I Math.} {\bf  325(7)} (1997) 749--752.


\bibitem{DLR} {\sc  B.\,Dehman, M.\,L\'eautaud and J.\,Le Rousseau}, Controllability of two coupled wave equations on a compact manifold. {\it Arch. Rational Mech. Anal.} {\bf 211(1)} (2014) 113--187.

\bibitem{EZ}
{\sc S.\,Ervedoza and E.\,Zuazua},  Sharp observability estimates for heat equations. {\it Arch. Rational Mech. Anal.}  {\bf 202(3)} (2011) 975--1017.



\bibitem{Ger}
{\sc P.\,G\' erard},  Microlocal Defect Measures. {\it Comm. Partial
Differential Equations}  {\bf 16(11)} (1991) 1761--1794.

\bibitem{Hor}
{\sc L.\,H\"ormander}, { On the uniqueness of the Cauchy problem. II}. {\it Math. Scand.} {\bf 7} (1959) 177--190.

\bibitem{Hor2}
{\sc L.\,H\" ormander}, {\it The Analysis of Linear Partial Differential Operators I-IV}, Springer, 1990


\bibitem{LZ}
{\sc M.\,Lazar and E.\,Zuazua},  Averaged control and observation of parameter-depending wave equations.  {\it C. R. Acad. Sci. Paris, Ser. I} {\bf 352(6)} (2014) 497--502.



\bibitem{LRTT}
{\sc G.\,Lebeau, J.\,Le Rousseau, P.\,Terpolilli and E.\,Tr\'elat}, { Some new results for the controllability of waves
equations}, presented on {\sl Workshop New trends in modeling, control and inverse problems}, Institut de Mathématiques de Toulouse, June 16 - 19, 2014. 
http://www.math.univ-toulouse.fr/$\sim$ervedoza/WebpageCIMI-Enrique/Slides/slides-lebeau.pdf

\bibitem{LRH}
{\sc T.\,Li, B.\,Rao and L.\,Hu}, Exact boundary synchronization for a coupled system of 1-D wave equations.   {\it ESAIM: COCV} {\bf 20} (2014) 339--361.


\bibitem{LuZ}
{\sc Q.\,L\"u and  E.\,Zuazua},  Averaged controllability for random evolution partial differential  equations, submitted.

\bibitem{Tar} {\sc L.\,Tartar},  H-measures, a new approach for studying homogenisation,
oscillation and concentration effects in PDEs. {\it Proc. Roy. Soc.
Edinburgh. Sect. A} {\bf 115(3-4)} (1990) 193--230.

\bibitem{Z-aver}
{\sc E.\,Zuazua},    Averaged Control.  {\it Automatica} {\bf 50(12)} (2014) 3077--3087.

\bibitem{Z}
{\sc E.\,Zuazua},   Robust observation of Partial Differential Equations, preprint (2014). 





\end{thebibliography}
\end{document}